\documentclass[numbers,webpdf]{ima-authoring-template}%

\graphicspath{{Fig/}}
\usepackage{placeins}
\usepackage[capitalize]{cleveref}

\theoremstyle{thmstyletwo}%
\newtheorem{theorem}{Theorem}%  meant for continuous numbers
\newtheorem{proposition}[theorem]{Proposition}%

\newtheorem{definition}{Definition}

\numberwithin{equation}{section}

\begin{document}

%\DOI{DOI HERE}
%\copyrightyear{2026}
%\vol{00}
%\pubyear{2026}
%\access{Advance Access Publication Date: Day Month Year}
%\appnotes{Paper}
%\copyrightstatement{Published by Oxford University Press on behalf of the Institute of Mathematics and its Applications. All rights reserved.}
%\firstpage{1}

%\subtitle{Subject Section}

\title[Residual-Guided Koopman Approximation]{Residual-Guided Dictionary Learning for\\Spectrally Accurate Koopman Approximation}

\author{George Coote*
\address{\orgdiv{Department of Applied Mathematics and Theoretical Physics}, \orgname{University of Cambridge}, \orgaddress{\street{Wilberforce Road}, \postcode{CB3 0WA}, \country{United Kingdom}}}}
\author{Matthew J. Colbrook\ORCID{0000-0003-4964-9575}
\address{\orgdiv{Department of Applied Mathematics and Theoretical Physics}, \orgname{University of Cambridge}, \orgaddress{\street{Wilberforce Road}, \postcode{CB3 0WA}, \country{United Kingdom}}}}
%\author{Fourth Author
%\address{\orgdiv{Department}, \orgname{Organization}, \orgaddress{\street{Street}, \postcode{Postcode}, \state{State}, \country{Country}}}}

\authormark{George Coote and Matthew Colbrook}

\corresp[*]{Corresponding author: \href{m.colbrook@damtp.cam.ac.uk}{m.colbrook@damtp.cam.ac.uk}}

%\received{Date}{0}{Year}
%\revised{Date}{0}{Year}
%\accepted{Date}{0}{Year}

%\editor{Associate Editor: Name}

\abstract{Koopman theory promises linear structure in nonlinear dynamics, but numerical Koopman spectra are easy to compute and hard to trust. A finite EDMD matrix always has eigenvalues; the problem is that many of them may have nothing to do with the infinite-dimensional operator. In this paper we make spectral reliability the objective of dictionary learning. We train neural-network dictionaries not merely to predict the next snapshot, but to minimize Residual Dynamic Mode Decomposition residuals: operator-level a posteriori errors that test whether computed eigenvalues and modes are genuine Koopman spectral objects. To keep the learned observables from collapsing into an unstable coordinate system, the loss also penalizes the condition number of the lifted data matrix. Thus the method couples two requirements that should not be separated: small Koopman residuals and a well-conditioned representation. The result is a learned dictionary that is expressive, numerically stable, and spectrally disciplined. Across conservative and dissipative benchmark systems, the method sharply reduces spectral pollution, improves residual pseudospectral inclusion, and lowers forecast error relative to standard fixed dictionaries. On sea-surface temperature data, it gives cleaner Koopman diagnostics and substantially better one-step forecasts from noisy observations with no governing equations. The message is simple: neural Koopman learning should be judged not by prediction alone, but by whether its spectral claims can be certified. Residuals provide the certificate; conditioning makes it computable.}
\keywords{Data-driven dynamics; Koopman operator; neural networks; Residual Dynamic Mode Decomposition.}

% \boxedtext{
% \begin{itemize}
% \item Key boxed text here.
% \item Key boxed text here.
% \item Key boxed text here.
% \end{itemize}}

\maketitle

\section{Introduction}

Koopman theory offers a striking promise: nonlinear dynamics studied through a linear operator on observables. If $F:X\to X$ advances the state of a dynamical system, the Koopman operator is the composition operator $\mathcal K_F g=g\circ F$. Introduced by Koopman and von Neumann in the 1930s \citep{koopman_1931,koopman_1932}, this viewpoint has been transformed by modern data-driven approximation methods \citep{mezic2005spectral,brunton2021modern}. It now appears across control, climate dynamics, epidemiology, neuroscience, non-autonomous systems, and interpretable machine learning \citep{haggerty2023control,froyland2021spectral,Orvieto2023RNN,Proctor2015DiseaseDMD,brunton2016,brunton2020machine}. The price of linearity is infinite dimension. Thus the numerical question is not how to compute a finite matrix, but whether that matrix has spectral meaning.

The question matters because Koopman spectra are used to identify coherent structures, frequencies, decay rates, metastable states, and reduced-order coordinates. Dynamic Mode Decomposition and Extended Dynamic Mode Decomposition approximate the Koopman operator by projection onto a finite dictionary of observables \cite{rowley2009,mezic2013,williams2015}. The dictionary is decisive. A bad choice may give accurate one-step predictions yet produce spurious eigenvalues, miss genuine spectral components, or return modes that are artifacts of the approximation.

There are two natural responses. One is to prescribe the dictionary: polynomials, Fourier modes, radial basis functions, Hermite functions, kernels, or delay coordinates. This can be powerful when the geometry is known, and brittle when it is not. The other is to learn the dictionary from data. Neural dictionaries have become a major direction in Koopman learning, from EDMD with dictionary learning \cite{boltt2017} and invariant-subspace and deep-DMD methods \cite{takeishi2017,yeung2019} to autoencoder and latent-linear models that seek coordinates with approximately linear dynamics \cite{lusch2018,otto2019,pan2020,logothetis2025,azencot2020}. Such methods add flexibility. They do not by themselves solve the spectral problem.

The obstruction is spectral pollution. A finite EDMD matrix always has eigenvalues, but they need not correspond to spectrum of $\mathcal K_F$. For infinite-dimensional operators, a convincing eigenvalue plot can be wrong, even if it is stable under different discretizations. What is needed is an a posteriori test: a computed spectral object should be trusted because it has a small residual for the Koopman operator itself, not merely because it is an eigenpair of a finite projection.

Residual Dynamic Mode Decomposition supplies such a test \cite{colbrook2023,colbrook2025}. ResDMD augments EDMD with residuals that converge, in the large-data limit, to residuals for the infinite-dimensional Koopman operator. It also computes residual pseudospectra: regions where the dictionary contains observables with small Koopman residual. Thus spectral pollution can be detected, and spectral inclusion assessed. The framework has now been developed for robust Koopman spectral computation, stochastic dynamics, snapshot-limited data, transfer operators, and scientific applications \cite{colbrookTownsend2023,colbrook2023,colbrookLiRautTownsend2023,colbrook2024fewSnapshots,schurig2026,colbrook2025,lorenzo2026residual,preprint2_to_appear}.

This paper joins these two threads. We learn neural dictionaries, but train them with ResDMD residuals rather than prediction error alone. The goal is not just a low-dimensional predictor, but a dictionary whose Koopman eigenvalues, eigenfunctions, residuals, and residual pseudospectra have numerical meaning. Learning improves the dictionary; residuals keep the spectral claims honest. Conditioning is crucial, but it has not yet been made part of the learning objective. Residual minimization alone can produce nearly linearly dependent observables. Then the lifted data matrix $\mathbf W^{1/2}\Psi_X$ is ill-conditioned, the Gram solve defining EDMD becomes unstable, and both eigenvalues and residuals may be fragile. Our loss couples the ResDMD residual with an explicit penalty on $\kappa(\mathbf W^{1/2}\Psi_X)$: the residual term seeks spectrally accurate observables, while the condition-number term prevents collapse to an unstable coordinate system. The result is a Koopman approximation that is expressive, stable, and computable.

The numerical effect is substantial. Across benchmarks, the trained dictionaries give smaller forecast errors, better-conditioned lifted data matrices, cleaner residual pseudospectral inclusion, and less spectral pollution. For conservative systems, the learned residual pseudospectra recover the expected unit-circle structure while suppressing spurious interior eigenvalues. For dissipative systems and real data, the same training principle improves both spectral diagnostics and prediction. A sea-surface-temperature example shows that these benefits persist for noisy observations, even when no governing equation is available.

The message is simple. In Koopman learning, neural dictionaries should be judged not just by expressivity or prediction error, but by spectral trustworthiness. ResDMD residuals supply the spectral target, and conditioning determines whether the computation is numerically sound. Optimizing both gives learned Koopman approximations that are accurate, stable, and informative beyond what fixed dictionaries alone can deliver.

\section{Mathematical Prerequisites}
\subsection{Koopman Operators}

We begin by recalling the definition of a Koopman operator; see, for example, \cite[Chapter~11]{Colbrook2026InfiniteDimensionalSpectralComputations}. Let $X\subseteq\mathbb R^d$ be the state space of a dynamical system, and let $x\in X$ denote its state at a given time, for example the position of a particle in a fluid. Let $\omega$ be a non-atomic, $\sigma$-finite Borel measure on $X$, not necessarily Lebesgue measure, and let $F:X\to X$ be Borel measurable.

The pushforward of $\omega$ by $F$ is the Borel measure
$$
    F_\#\omega(S)=\omega(F^{-1}(S)), \qquad S\subseteq X \text{ Borel}.
$$
We say that $F$ is non-singular with respect to $\omega$ if
$$
    \omega(S)=0 \quad \Longrightarrow \quad F_\#\omega(S)=0,
$$
or equivalently if $F_\#\omega\ll\omega$. We use the following definition.

\begin{definition}
Let $X$ and $\omega$ be as above, and let $F:X\to X$ be a non-singular Borel map. The Koopman operator associated with $F$ is the operator $\mathcal K_F$ on $L^\infty(X,\omega)$ defined by
$$
    (\mathcal K_F g)(x)=g(F(x)), \qquad x\in X.
$$
Here $g$ denotes any representative of its equivalence class. The non-singularity of $F$ makes this definition independent of the representative.
\end{definition}

By the Radon--Nikodym theorem, there is a Borel measurable density $\rho_F:X\to[0,\infty]$ such that
$$
    F_\#\omega(S)=\int_S \rho_F\,d\omega,\qquad S\subseteq X \text{ Borel}.
$$
Thus $F_\#\omega$ is determined by $\rho_F$ relative to $\omega$. The operator $\mathcal K_F$ is bounded on $L^\infty(X,\omega)$. For $1\leq p<\infty$, however, this is no longer automatic: $\mathcal K_F$ need not even map $L^p(X,\omega)$ into itself.

\begin{proposition}
Suppose that $\rho_F\in L^\infty(X,\omega)$. Then, for every $1\leq p<\infty$, the Koopman operator $\mathcal K_F$ is bounded on $L^p(X,\omega)$, and
$
    \|\mathcal K_F\|_{L^p\to L^p}=\|\rho_F\|_\infty^{1/p}.
$
\end{proposition}

Throughout the paper we assume that $\mathcal K_F$ is a bounded linear operator on the Hilbert space $L^2(X,\omega)$. We say that $F:X\to X$ is measure preserving if $F_\#\omega=\omega$. If, in addition, $\mathcal K_F$ is invertible on $L^2(X,\omega)$, we call the associated system invertible. For a measure-preserving map $F$, the Koopman operator is an isometry on $L^2(X,\omega)$. Thus, if the system is invertible, $\mathcal K_F$ is unitary, and hence $\operatorname{Sp}(\mathcal K_F)\subseteq\mathbb T$.

\subsection{Dictionaries}

Since $L^2(X,\omega)$ is infinite-dimensional, $\mathcal K_F$ is not represented by a finite matrix. Numerical computation therefore begins with a finite-dimensional surrogate. We choose a finite set of observables $\mathcal D_N=(\psi_1,\ldots,\psi_N)$ and approximate the compression of $\mathcal K_F$ onto $V_N=\operatorname{span}\{\psi_1,\ldots,\psi_N\}$. We call $\mathcal D_N$ a dictionary. No orthogonality or normalization is assumed. The associated feature map is
$$
    \Psi(x)=(\psi_1(x),\ldots,\psi_N(x)),
$$
which lifts the state $x\in X$ to a point in the feature space $\mathbb F^N$, where $\mathbb F=\mathbb R$ or $\mathbb C$.

We judge finite-dimensional approximations of $\mathcal K_F$ by their spectral accuracy. For this purpose, the dictionary should not only approximate typical observables, but also resolve eigenfunctions and pseudoeigenfunctions. This choice usually reflects information about the dynamics. A periodic dictionary, for example, is natural for highly oscillatory eigenfunctions, but may give poor approximations to eigenfunctions with rapid decay.

We first give several one-dimensional examples, then describe how to combine them into dictionaries in higher dimensions.

\subsubsection{Fourier Modes}
A sufficiently regular $1$-periodic function has the Fourier expansion
$$
    f(x)=a_0+\sum_{n=1}^{\infty}\bigl(a_n\cos(2\pi n x)+b_n\sin(2\pi n x)\bigr).
$$
Small values of $n$ describe slowly varying structure, while large values describe high-frequency oscillations. Thus, for a reasonably smooth periodic function, a modest number of Fourier modes may already give an accurate approximation. This motivates the one-dimensional trigonometric dictionary formed from $1$, $\cos(2\pi n x)$, and $\sin(2\pi n x)$.

\subsubsection{Chebyshev Polynomials}
For nonperiodic problems on a finite interval, a standard analogue of Fourier modes is the Chebyshev polynomial basis. The Chebyshev polynomial $T_n$ has degree $n$ and is characterized by
$$
    T_n(\cos\theta)=\cos(n\theta), \qquad \theta\in\mathbb R.
$$
The polynomials may be evaluated by the three-term recurrence
$$
    T_{n+1}(x)=2xT_n(x)-T_{n-1}(x), \qquad n\geq 1,
$$
with $T_0(x)=1$ and $T_1(x)=x$. We use $T_0,T_1,T_2,\ldots$ as one-dimensional dictionary functions, after rescaling the physical interval to $[-1,1]$ if necessary.

\subsubsection{Hermite Functions}
On the real line, a natural analogue is provided by the Hermite functions. They are the eigenfunctions of the quantum harmonic oscillator and are closely tied to Gauss--Hermite quadrature. Let $P_n$ denote the physicists' Hermite polynomials, defined by
$$
    P_{n+1}(x)=2xP_n(x)-2nP_{n-1}(x), \qquad n\geq 1,
$$
with $P_0(x)=1$ and $P_1(x)=2x$. The corresponding Hermite functions are, up to normalization,
$$
    h_n(x)=P_n(x)e^{-x^2/2}, \qquad x\in\mathbb R.
$$
We use the functions $h_0,h_1,h_2,\ldots$ as one-dimensional dictionary functions.

\subsubsection{Forming Higher-Dimensional Dictionaries}
Let $f,g\in L^2(\mathbb R)$, or more generally $f,g\in L^2(I)$ for an interval $I$. Their tensor product is
$
    (f\otimes g)(x,y)=f(x)g(y).
$
Thus $f\otimes g\in L^2(\mathbb R^2)$, with the obvious modification on product intervals. If $\mathcal B_1$ and $\mathcal B_2$ are orthonormal bases for $L^2(\mathbb R^{d_1})$ and $L^2(\mathbb R^{d_2})$, respectively, then
$
    \{f\otimes g:f\in\mathcal B_1,\ g\in\mathcal B_2\}
$
is an orthonormal basis for $L^2(\mathbb R^{d_1+d_2})$. Analogous tensor-product constructions apply in many other approximation spaces.

The same idea applies to finite dictionaries. One may tensor dictionaries of the same type in each coordinate, or choose different one-dimensional dictionaries to reflect different variables. For example, in a pendulum problem the state consists of an angle and a velocity. A trigonometric dictionary is natural for the angular variable, while a Hermite-type dictionary is more appropriate for the unbounded velocity variable.

\subsection{ResDMD}
We summarize Residual Dynamic Mode Decomposition (ResDMD) \cite{colbrook2023}, following \cite{colbrook2025}. ResDMD augments Extended Dynamic Mode Decomposition (EDMD): in addition to approximate eigenpairs, it computes residuals. These residuals quantify the error in the computed eigenpairs and, in the large-data limit, converge to the corresponding infinite-dimensional residuals.

ResDMD is built from snapshot pairs. These may come from one long trajectory or from many shorter ones. Here we sample a collection of initial conditions at random and evolve each for a prescribed number of time steps. The data are written as
\begin{equation}\label{snapshot}
    \bigl(x^{(m)},y^{(m)}\bigr)_{m=1}^M,\qquad y^{(m)}=F\bigl(x^{(m)}\bigr),
\end{equation}
where $y^{(m)}$ is the image of $x^{(m)}$ after one time step.

Fix a Koopman operator $\mathcal K=\mathcal K_F$, and suppose that we are given the snapshot pairs in \cref{snapshot}. Let $\mathcal D=(\psi_j)_{j\geq 1}$ be a linearly independent dictionary in $L^2(X,\omega)$. For $N\geq 1$, set $V_N=\operatorname{span}\{\psi_1,\ldots,\psi_N\}$, and let $\mathcal P_N:L^2(X,\omega)\to V_N$ denote the orthogonal projection. We seek the matrix $\mathbb K_N\in\mathbb C^{N\times N}$ representing the compression $\mathcal P_N\mathcal K|_{V_N}$ in this basis. Thus, for $i=1,\ldots,N$,
$$
    \mathcal K\psi_i=\sum_{j=1}^N [\mathbb K_N]_{ji}\psi_j+\rho_i,
$$
where $\rho_i\perp V_N$ is the projection residual. Equivalently, the $i$th column of $\mathbb K_N$ minimizes
$$
    \int_X \left|\sum_{j=1}^N [\mathbb K_N]_{ji}\psi_j(x)-\psi_i(F(x))\right|^2\,d\omega(x).
$$
The associated normal equations are
$$
    \sum_{j=1}^N \langle \psi_j,\psi_\ell\rangle [\mathbb K_N]_{ji}=\langle \mathcal K\psi_i,\psi_\ell\rangle,\qquad \ell=1,\ldots,N.
$$
Writing $[G_N]_{\ell j}=\langle \psi_j,\psi_\ell\rangle$, the linear independence of the dictionary functions implies that $G_N$ is Hermitian positive definite, and hence $\mathbb K_N$ is uniquely determined.

In the data-driven setting the inner products above are not known exactly. Instead, we approximate the relevant integrals by a quadrature rule on the snapshot points,
$$
    \int_X |\rho_i(x)|^2\,d\omega(x)\approx \sum_{m=1}^M w_m |\rho_i(x^{(m)})|^2,
$$
with weights $w_m>0$. We assume that this quadrature converges in the large-data limit.

Let $\mathbf W=\operatorname{diag}(w_1,\ldots,w_M)$, and recall the feature map $\Psi(x)=(\psi_1(x),\ldots,\psi_N(x))$. The lifted data matrices are
$$
    \Psi_X=
    \begin{pmatrix}
        \Psi(x^{(1)})\\
        \vdots\\
        \Psi(x^{(M)})
    \end{pmatrix},
    \qquad
    \Psi_Y=
    \begin{pmatrix}
        \Psi(y^{(1)})\\
        \vdots\\
        \Psi(y^{(M)})
    \end{pmatrix}
    \in\mathbb C^{M\times N}.
$$
Thus $[\Psi_X]_{mj}=\psi_j(x^{(m)})$ and $[\Psi_Y]_{mj}=\psi_j(y^{(m)})$. EDMD chooses $\mathbf K=\mathbf K_N\in\mathbb C^{N\times N}$ to solve the weighted least-squares problem
$$
    \min_{\mathbf K\in\mathbb C^{N\times N}}\sum_{m=1}^M w_m\left\|\Psi(y^{(m)})-\Psi(x^{(m)})\mathbf K\right\|_{\ell^2}^2
    =
    \min_{\mathbf K\in\mathbb C^{N\times N}}\left\|\mathbf W^{1/2}\Psi_Y-\mathbf W^{1/2}\Psi_X\mathbf K\right\|_F^2.
$$
The Moore--Penrose solution is
$
    \mathbf K=(\mathbf W^{1/2}\Psi_X)^\dagger(\mathbf W^{1/2}\Psi_Y).
$
If $\mathbf W^{1/2}\Psi_X$ has full column rank, this reduces to
$
    \mathbf K=(\Psi_X^\ast \mathbf W\Psi_X)^{-1}\Psi_X^\ast \mathbf W\Psi_Y.
$
In applications, $\mathbf W^{1/2}\Psi_X$ is often ill-conditioned. We therefore use the Tikhonov-regularized EDMD matrix
$$
    \mathbf K_\varepsilon=\bigl(\Psi_X^\ast \mathbf W\Psi_X+\varepsilon\mathbf I\bigr)^{-1}\Psi_X^\ast \mathbf W\Psi_Y,
$$
with $\varepsilon>0$. The ridge parameter should be large enough to stabilize the inversion of the empirical Gram matrix $\Psi_X^\ast \mathbf W\Psi_X$, but small enough that it does not dominate the least-squares fit.

For a fixed, untrained dictionary, ill-conditioning usually has one of two causes: the dictionary functions are nearly linearly dependent on the sampled data, or the sampled trajectories are themselves highly correlated. The first issue can be addressed by pruning the dictionary. The second can often be reduced by using a more diverse set of trajectories, though for highly regular systems some correlation is unavoidable.

The eigenpairs of $\mathbf K$ provide candidate Koopman eigenpairs. To assess them, we use residuals. If $\mathcal K$ is normal, then for every $z\in\mathbb C$,
$$
    \operatorname{dist}\bigl(z,\operatorname{Sp}(\mathcal K)\bigr)=\inf_{g\ne 0}\frac{\|(\mathcal K-zI)g\|}{\|g\|}.
$$
Thus a small residual is direct evidence that $z$ is close to the spectrum.

Let $g\in V_N$ have the expansion
$$
    g=\sum_{i=1}^N g_i\psi_i,\qquad \mathbf g=(g_1,\ldots,g_N)^\top.
$$
The empirical residual associated with $z\in\mathbb C$ and $\mathbf g\ne 0$ is
\begin{equation}
\label{eqn:residual}
    \operatorname{res}(z,\mathbf g)=
    \frac{\|(\mathbf W^{1/2}\Psi_Y-z\mathbf W^{1/2}\Psi_X)\mathbf g\|_{\ell^2}}{\|\mathbf W^{1/2}\Psi_X\mathbf g\|_{\ell^2}}.
\end{equation}
For numerical optimization it is often preferable to work with $\operatorname{res}(z,\mathbf g)^2$, avoiding the square root. Under the assumed convergence of the quadrature rule, and for fixed $g\in V_N$,
$$
    \operatorname{res}(z,\mathbf g)\longrightarrow \frac{\|(\mathcal K-zI)g\|}{\|g\|}\qquad\text{as }M\to\infty.
$$
We therefore use
$
    \tau(z)=\inf_{\mathbf g\ne 0}\operatorname{res}(z,\mathbf g)
$
as a residual pseudospectral diagnostic. This interpretation gives a direct distance-to-spectrum certificate when $\mathcal K$ is normal; for non-normal operators it remains useful, but should be read with more caution.

It remains to compute $\tau(z)$. Assume that $\mathbf W^{1/2}\Psi_X$ has full column rank, and take its economy QR factorization
$
    \mathbf W^{1/2}\Psi_X=\mathbf Q\mathbf R,
$
where $\mathbf Q^\ast\mathbf Q=\mathbf I$ and $\mathbf R$ is upper triangular with positive diagonal entries. Setting $\mathbf w=\mathbf R\mathbf g$, we have
$
    \|\mathbf W^{1/2}\Psi_X\mathbf g\|_{\ell^2}^2=\|\mathbf w\|_{\ell^2}^2.
$
Hence
$$
    \operatorname{res}(z,\mathbf g)^2=
    \frac{\|(\mathbf W^{1/2}\Psi_Y\mathbf R^{-1}-z\mathbf Q)\mathbf w\|_{\ell^2}^2}{\|\mathbf w\|_{\ell^2}^2}.
$$
Taking the infimum over $\mathbf g\ne 0$, equivalently over $\mathbf w\ne 0$, gives
$$
    \tau(z)=\sigma_{\min}\bigl(\mathbf W^{1/2}\Psi_Y\mathbf R^{-1}-z\mathbf Q\bigr)=\sqrt{\lambda_{\min}(M(z))},
$$
where
\begin{equation}
\label{eqn:mz}
    M(z)=(\mathbf R^\ast)^{-1}\Psi_Y^\ast \mathbf W\Psi_Y\mathbf R^{-1}
    -z(\mathbf R^\ast)^{-1}\Psi_Y^\ast \mathbf W^{1/2}\mathbf Q
    -\overline z\,\mathbf Q^\ast\mathbf W^{1/2}\Psi_Y\mathbf R^{-1}
    +|z|^2\mathbf I.
\end{equation}
This is the Hermitian positive semidefinite matrix obtained by expanding
$$
    \bigl(\mathbf W^{1/2}\Psi_Y\mathbf R^{-1}-z\mathbf Q\bigr)^\ast
    \bigl(\mathbf W^{1/2}\Psi_Y\mathbf R^{-1}-z\mathbf Q\bigr),
$$
using $\mathbf Q^\ast\mathbf Q=\mathbf I$. \cref{algorithm:resdmd,algorithm:pseudospectra} summarize the computational procedures.

\renewcommand{\algorithmicrequire}{\textbf{Input:}}
\renewcommand{\algorithmicensure}{\textbf{Output:}}

\begin{algorithm}[t]
\caption{ResDMD for computing eigendecomposition and residuals}\label{algorithm:resdmd}
\begin{algorithmic}[1]
\Require Snapshot data $((x^{(m)},y^{(m)}))_{m=1}^M$, quadrature weights $(w_m)_{m=1}^M$, dictionary $(\psi_j)_{j=1}^N$ with $M>N$, and ridge parameter $\varepsilon>0$.
\State Compute the lifted data matrices $\Psi_X$, $\Psi_Y$, and set $\mathbf W=\operatorname{diag}(w_1,\ldots,w_M)$.
\State Compute $\mathbf K_\varepsilon=(\Psi_X^\ast\mathbf W\Psi_X+\varepsilon\mathbf I)^{-1}\Psi_X^\ast\mathbf W\Psi_Y$.
\State Compute the eigendecomposition $\mathbf K_\varepsilon\mathbf V=\mathbf V\Lambda$.
\State For each eigenpair $(\lambda_j,\mathbf v_j)$, compute
$$
    r_j=\operatorname{res}(\lambda_j,\mathbf v_j)=\frac{\|(\mathbf W^{1/2}\Psi_Y-\lambda_j\mathbf W^{1/2}\Psi_X)\mathbf v_j\|_{\ell^2}}{\|\mathbf W^{1/2}\Psi_X\mathbf v_j\|_{\ell^2}}.
$$
\Ensure Eigenvalues $\Lambda$, eigenvectors $\mathbf V\in\mathbb C^{N\times N}$, and residuals associated with $x\mapsto\Psi(x)\mathbf v_j$.
\end{algorithmic}
\end{algorithm}

\begin{algorithm}[t]
\caption{ResDMD for computing pseudospectra}\label{algorithm:pseudospectra}
\begin{algorithmic}[1]
\Require Snapshot data $((x^{(m)},y^{(m)}))_{m=1}^M$, quadrature weights $(w_m)_{m=1}^M$, dictionary $(\psi_j)_{j=1}^N$ with $M>N$, and grid points $(z_\ell)_{\ell=1}^k$.
\State Compute the lifted data matrices $\Psi_X$, $\Psi_Y$, and set $\mathbf W=\operatorname{diag}(w_1,\ldots,w_M)$.
\State Compute an economy QR factorization $\mathbf W^{1/2}\Psi_X=\mathbf Q\mathbf R$, where $\mathbf Q\in\mathbb C^{M\times N}$ and $\mathbf R\in\mathbb C^{N\times N}$.
\State For each grid point $z_\ell$, compute $M(z_\ell)$ defined in \eqref{eqn:mz}.
\State Compute $\tau(z_\ell)=\sqrt{\lambda_{\min}(M(z_\ell))}$ for each grid point $z_\ell$.
\Ensure The values $\{\tau(z_\ell):1\leq\ell\leq k\}$.
\end{algorithmic}
\end{algorithm}

\subsection{Neural Networks}
Neural networks provide a flexible way to approximate eigenfunctions. Once trained, they are cheap to evaluate, have a simple compositional structure, and possess strong approximation properties; see, for example, \cite[pp.~208--209]{haykin1999}. We recall the basic definition.

\begin{definition}
Let $X\subseteq\mathbb R^d$. A feedforward neural network is a map $f_\theta:X\to\mathbb R^q$ of the form
$$
    f_\theta=\Phi_L\circ\Phi_{L-1}\circ\cdots\circ\Phi_1\big|_X,
$$
where each layer map $\Phi_\ell$ is determined by trainable weights and biases. Typically,
$$
    \Phi_\ell(z)=\sigma(A_\ell z+b_\ell)
$$
for the hidden layers, with an affine or activated affine map at the output.
\end{definition}

The maps $\Phi_\ell$ are called layers. The first layer receives the $d$-dimensional input, the last layer produces the output, and the intermediate layers are hidden layers. The common dimension of the hidden layers is called the hidden dimension, or width. The activation function, together with the number and widths of the layers, specifies the architecture of the network.

Throughout the paper we use the rectified linear unit activation,
$$
    \operatorname{ReLU}(x)=
    \begin{cases}
        x, & x\geq 0,\\
        0, & x<0.
    \end{cases}
$$
For a vector $\mathbf x=(x_1,\ldots,x_n)\in\mathbb R^n$, ReLU is applied componentwise:
$$
    \operatorname{ReLU}(\mathbf x)=\bigl(\operatorname{ReLU}(x_1),\ldots,\operatorname{ReLU}(x_n)\bigr).
$$
The appeal of ReLU is twofold: it is cheap to evaluate, and it tends to produce sparse activations, with many hidden units returning zero \cite{glorot2011}. In our benchmarks, ReLU layers performed particularly well.

We now specify the architecture used in the computations. Fix an input dimension $d$, a hidden dimension $d_h$, and an integer $n\geq 1$. The network has layer maps $\varsigma_0,\varsigma_1,\ldots,\varsigma_n$ and scalar output. The first layer is
$$
    \varsigma_0(x)=\operatorname{ReLU}(W_0^\top x+b_0),\qquad W_0\in\mathbb R^{d\times d_h},\quad b_0\in\mathbb R^{d_h}.
$$
For $1\leq i\leq n-1$, the hidden layers are
$$
    \varsigma_i(x)=\operatorname{ReLU}(W_i^\top x+b_i),\qquad W_i\in\mathbb R^{d_h\times d_h},\quad b_i\in\mathbb R^{d_h}.
$$
The final layer is affine,
$$
    \varsigma_n(x)=W_n^\top x+b_n,\qquad W_n\in\mathbb R^{d_h},\quad b_n\in\mathbb R.
$$
No activation is applied at the output, so the network is not artificially restricted in sign or magnitude. The resulting network is
$$
    f(x)=(\varsigma_n\circ\varsigma_{n-1}\circ\cdots\circ\varsigma_0)(x),\qquad x\in X.
$$
Although this notation allows many architectures, the numerical experiments below use shallow networks with widths in the range $8\leq d_h\leq 32$.

To begin the optimization, the weights $W_i$ and biases $b_i$ must be initialized. We use the default PyTorch initialization for affine layers: entries of $W_i$ and $b_i$ are drawn independently from
$
    U(-\tfrac{1}{\sqrt{k_i}},\tfrac{1}{\sqrt{k_i}}),
$
where $k_i$ is the input dimension of the layer $\varsigma_i$. Other initializations are possible. A common alternative for ReLU networks is the Kaiming--He normal initialization,
$
    N(0,\tfrac{2}{k_i}).
$

\section{Setup}

The choice of dictionary usually requires prior information about the system. Common fixed choices include radial basis functions, Chebyshev polynomials, and Hermite functions. We call such choices \emph{untrained dictionaries}: they contain no trainable parameters.

Even when approximations based on these dictionaries converge, their spectra may exhibit severe pollution, producing numerical modes with no counterpart in the true dynamics. Such eigenvalues must then be removed, for example by residual filtering. At the same time, fixed dictionaries often give good trajectory predictions. This suggests a useful distinction: they may capture local prediction well while introducing unreliable global spectral features.

We instead consider \emph{trained dictionaries}: dictionaries whose elements depend on trainable parameters. Their appeal is that little prior knowledge of the dynamics is required, beyond the choice of architecture and a few hyperparameters. This makes them natural for complex systems whose structure is not known in advance.

\subsection{The Loss Function}
\label{sub:loss}

Let $N\in\mathbb N$, and let $\mathcal D_\Theta=(\psi_1,\ldots,\psi_N)$ be a dictionary of neural networks with a common architecture. The networks are initialized independently from the same distribution. We write $\theta_j$ for the parameter vector of $\psi_j$, and collect the parameters as $\Theta=(\theta_1,\ldots,\theta_N)$. The value of the $j$th network at $x$ is denoted by $\psi_j(x;\theta_j)$, and the corresponding feature map is
$$
    \Psi(x;\Theta)=\bigl(\psi_1(x;\theta_1),\ldots,\psi_N(x;\theta_N)\bigr).
$$
We assume snapshot data
$$
    \mathbf X=\bigl(x^{(m)},y^{(m)}\bigr)_{m=1}^M,\qquad N\ll M,
$$
obtained from simulation or experiment, together with quadrature weights $(w_m)_{m=1}^M$. Set $\mathbf W=\operatorname{diag}(w_1,\ldots,w_M)$. For a fixed parameter vector $\Theta$, define the lifted data matrices
$$
    \Psi_X(\Theta)=
    \begin{pmatrix}
        \Psi(x^{(1)};\Theta)\\
        \vdots\\
        \Psi(x^{(M)};\Theta)
    \end{pmatrix},
    \qquad
    \Psi_Y(\Theta)=
    \begin{pmatrix}
        \Psi(y^{(1)};\Theta)\\
        \vdots\\
        \Psi(y^{(M)};\Theta)
    \end{pmatrix}
    \in\mathbb C^{M\times N}.
$$

Applying ResDMD, as in \cref{algorithm:resdmd}, to the data $(\mathbf X,\mathbf W)$ with lifted matrices $\Psi_X(\Theta)$ and $\Psi_Y(\Theta)$ gives an EDMD matrix $\mathbf K(\Theta)\in\mathbb C^{N\times N}$, eigenpairs $\bigl(\lambda_j(\Theta),\mathbf v_j(\Theta)\bigr)_{j=1}^N$, and residuals
$
    r_j(\Theta)=\operatorname{res}\bigl(\lambda_j(\Theta),\mathbf v_j(\Theta)\bigr).
$
Equivalently, $\mathbf v_j(\Theta)$ represents the observable $x\mapsto \Psi(x;\Theta)\mathbf v_j(\Theta)$. We also define
$
    \kappa(\Theta)=\kappa\bigl(\mathbf W^{1/2}\Psi_X(\Theta)\bigr),
$
the condition number of the weighted lifted data matrix. Thus the EDMD matrix, its eigenpairs, the residuals, and the conditioning all depend on $\Theta$ through the trained dictionary.

We define the loss
\begin{equation}
\label{eqn:loss}
    L(\mathbf X;\Theta)=\frac{1}{N}\sum_{j=1}^N r_j(\Theta)^2+\alpha\,\operatorname{ReLU}\left(\log\frac{\kappa(\Theta)}{\kappa_0}\right),
\end{equation}
where $\alpha>0$ and $\kappa_0>1$ are hyperparameters. The first term penalizes large spectral residuals. The second penalizes ill-conditioning, but only once $\kappa(\Theta)$ exceeds the prescribed threshold $\kappa_0$.

Training approximately minimizes $L(\mathbf X;\Theta)$ using the AdamW optimizer described in \cref{sub:training}:
$$
    \Theta^\star\approx\operatorname*{arg\,min}_{\Theta}L(\mathbf X;\Theta).
$$

The first term in \eqref{eqn:loss} measures spectral accuracy. Ideally, one would minimize the squared distance to $\operatorname{Sp}(\mathcal K)$, but this quantity is not available. For normal $\mathcal K$, the operator residual gives an upper bound on this distance, and the ResDMD residuals converge, in the large-data limit, to the corresponding infinite-dimensional residuals. Thus a small mean-squared residual indicates a more reliable spectral approximation.

The second term penalizes ill-conditioning of $\mathbf W^{1/2}\Psi_X(\Theta)$. A large condition number indicates that the learned dictionary functions are nearly linearly dependent on the sampled data. Equivalently, the empirical Gram matrix $\Psi_X(\Theta)^\ast\mathbf W\Psi_X(\Theta)$ is close to singular, so the least-squares problem defining the EDMD matrix is ill-conditioned. The penalty therefore serves two purposes: it improves numerical stability and encourages the dictionary to span a full-rank empirical subspace.

We use a logarithmic penalty rather than a linear or quadratic one. Since $\kappa$ is a ratio, $\log\kappa$ places conditioning on an additive scale and avoids the excessive growth of polynomial penalties for large $\kappa$. The ReLU gate makes the penalty inactive below the threshold $\kappa_0$. Thus the loss discourages serious ill-conditioning without penalizing harmless variation in the condition number.

\begin{algorithm}[t!]
\caption{Training loop}
\label{algorithm:training}
\begin{algorithmic}[1]
\Require Training data $\mathbf X_0$, test data $\mathbf X_1$, both with weights, initial parameters $\Theta$, and hyperparameters $(\alpha,\kappa_0,\eta,p,E_{\max},P,\gamma,\varepsilon)$.
\Ensure Trained parameters $\Theta^\star$.
\State $L^\star\gets+\infty$, $\Theta^\star\gets\Theta$, $c\gets 0$.
\For{$e=1,\ldots,E_{\max}$}
    \State Sample $\mathbf X_0'\subset\mathbf X_0$ of size $\lfloor p|\mathbf X_0|\rfloor$ without replacement.
    \State Compute $\Psi_X(\Theta)$ and $\Psi_Y(\Theta)$ on $\mathbf X_0'$.
    \State Compute $\mathbf K(\Theta)$ via \cref{algorithm:resdmd}.
    \State Compute eigenpairs $(\lambda_j(\Theta),\mathbf v_j(\Theta))$, residuals $r_j(\Theta)$, and condition number $\kappa(\Theta)$.
    \State Evaluate $L(\mathbf X_0';\Theta)$ via \eqref{eqn:loss}.
    \State Compute $\nabla_\Theta L(\mathbf X_0';\Theta)$ using autodifferentiation.
    \State Update $\Theta\gets\mathrm{AdamW}(\Theta,\nabla_\Theta L,\eta)$.
    \State Evaluate test loss $L(\mathbf X_1;\Theta)$.
    \If{$L(\mathbf X_1;\Theta)<(1-\gamma)L^\star$}
        \State $L^\star\gets L(\mathbf X_1;\Theta)$, $\Theta^\star\gets\Theta$, $c\gets 0$.
    \Else
        \State $c\gets c+1$.
        \If{$c\geq P$}
            \State \textbf{break}
        \EndIf
    \EndIf
\EndFor
\State \Return $\Theta^\star$.
\end{algorithmic}
\end{algorithm}

\subsection{Training Process}
\label{sub:training}

We minimize $L(\mathbf X;\Theta)$, defined in \cref{sub:loss}, by mini-batch stochastic optimization with AdamW. The data consist of snapshot triples $(x^{(m)},y^{(m)},w_m)_{m=1}^M$, where $w_m$ is the quadrature weight associated with the pair $(x^{(m)},y^{(m)})$. We shuffle these triples uniformly at random and split them into a training set $\mathbf X_0$, containing $70\%$ of the snapshots, and a validation set $\mathbf X_1$, containing the remaining $30\%$. This keeps each snapshot pair attached to its weight.

After the split, the weighted sums are used as empirical training and test objectives, rather than as independent quadrature rules. The held-out set is used to monitor generalization and to implement early stopping. Unless stated otherwise, we use the following hyperparameters:
\begin{enumerate}
    \item condition-number penalty weight $\alpha=10^{-2}$ and threshold $\kappa_0=10^3$;
    \item AdamW learning rate $\eta=10^{-3}$;
    \item mini-batch fraction $p=0.1$ of $\mathbf X_0$ per optimization step;
    \item maximum number of epochs $E_{\max}=500$;
    \item early stopping after $P=30$ epochs with no improvement in test loss, where improvement means a relative decrease of at least $r=0$;
    \item candidate ridge parameters $\varepsilon\in\{0.3,3,30,300\}$.
\end{enumerate}

\section{Results for Benchmark Systems}

We now present results for three benchmark systems. \cref{tab:benchmark_summary} summarizes the systems and the main numerical gains obtained from the trained dictionaries.

\begin{table}[t]
\centering
\small
\setlength{\tabcolsep}{6pt}
\caption{Summary of benchmark systems (rows 1 to 3) and real-world climate data (row 4). The condition number and one-step forecast error are reported for the trained and benchmark dictionaries; smaller values are better. The symbol $\pm$ denotes one standard deviation over the runs.}
\label{tab:benchmark_summary}
\begin{tabular}{p{0.18\textwidth}p{0.17\textwidth}p{0.25\textwidth}p{0.25\textwidth}}
\hline
System & Benchmark dictionary & Condition number & One-step forecast error \\
\hline
Pendulum & Fourier--Hermite &
\begin{tabular}[t]{@{}l@{}}trained: $(9.06 \pm 0.81)\times10^2$\\benchmark: $(1.16 \pm 0.12) \times 10^1$\end{tabular} &
\begin{tabular}[t]{@{}l@{}}trained: $(5.03 \pm 0.57) \times 10^{-4}$\\benchmark: $(6.18 \pm 0.54) \times 10^{-1}$\end{tabular} \\
Undamped oscillator & Hermite &
\begin{tabular}[t]{@{}l@{}}trained: $(2.59 \pm 0.40) \times 10^2$\\benchmark: $(4.45 \pm 0.60) \times 10^6$\end{tabular} &
\begin{tabular}[t]{@{}l@{}}trained: $(8.47 \pm 1.09) \times 10^{-5}$\\benchmark: $(4.91 \pm 0.62) \times 10^{-3}$\end{tabular} \\
Duffing oscillator & Chebyshev &
\begin{tabular}[t]{@{}l@{}}trained: $(7.17 \pm 0.21) \times 10^2$\\benchmark: $(4.75 \pm 0.72) \times 10^5$\end{tabular} &
\begin{tabular}[t]{@{}l@{}}trained: $(1.24 \pm 0.12) \times 10^{-4}$\\benchmark: $(2.10 \pm 0.20) \times 10^{-3}$\end{tabular} \\
Sea-Surface Temperature & Quadratic & 
\begin{tabular}[t]{@{}l@{}}trained ($5$ dim): $5.73 \pm 1.4$ \\ trained ($10$ dim): $5.48 \pm 2.03$ \\ benchmark ($5$ dim): $7.32$\\ benchmark ($10$ dim): $13.5$\end{tabular} &
\begin{tabular}[t]{@{}l@{}}trained ($5$ dim): $0.0433$\\ trained ($10$ dim): $0.0569 \pm 0.0277$ \\ benchmark ($5$ dim): $0.2194$\\ benchmark ($10$ dim): $0.3529$\end{tabular} \\
\hline
\end{tabular}
\end{table}

We plot the training and test losses over the training epochs. When the benchmark loss lies on a comparable scale, it is shown as a horizontal reference line. We give the analogous plot for the condition number. Test condition numbers are reported over $10$ runs, with $\pm$ denoting one standard deviation.

For each system we also show three ResDMD pseudospectral plots: one for a benchmark dictionary chosen using prior knowledge of the system, one for the randomly initialized neural-network dictionary before training, and one for the trained neural-network dictionary. Here ``pseudospectrum'' refers to the residual pseudospectrum computed in \cref{algorithm:pseudospectra}, rather than to the pseudospectrum of the finite EDMD matrix alone. At a grid point $z$, ResDMD computes
$$
    \tau(z)=\inf_{\mathbf g\ne 0}\frac{\|(\mathbf W^{1/2}\Psi_Y-z\mathbf W^{1/2}\Psi_X)\mathbf g\|_{\ell^2}}{\|\mathbf W^{1/2}\Psi_X\mathbf g\|_{\ell^2}}.
$$
Thus the $\delta$-level set
$$
    \Sigma_\delta^{M,N}=\{z\in\mathbb C:\tau(z)\leq\delta\}
$$
contains spectral parameters for which the dictionary spans an observable with empirical Koopman residual at most $\delta$. In the large-data limit this gives the corresponding infinite-dimensional residual test. In particular, if $\tau_\infty(z)\leq\delta$, then $z$ lies in the $\delta$-pseudospectrum of $\mathcal K$; when $\mathcal K$ is normal this gives the concrete inclusion
$$
    \operatorname{dist}\bigl(z,\operatorname{Sp}(\mathcal K)\bigr)\leq\delta.
$$
Small values of $\tau(z)$, shown in yellow, therefore indicate low-residual spectral candidates. EDMD eigenvalues are marked by red crosses. We interpret red crosses lying in large-$\tau(z)$ regions as spectral pollution, while low-$\tau(z)$ regions without nearby EDMD eigenvalues indicate poor spectral inclusion by the finite matrix. The pseudospectral figures show representative runs; scalar diagnostics are reported across multiple runs where stated.

Finally, we report the normalized one-step forecast error $\|\Psi_Y-\Psi_X\mathbf K\|_F^2/\|\Psi_Y\|_F^2$, again over $15$ runs. The normalization is important because different dictionary families may operate on very different scales. For example, Hermite dictionaries can take much larger values than the trained neural-network dictionaries, so absolute errors alone are not directly comparable. The same convention for standard deviations is used for the forecast errors.

All computations were performed on an NVIDIA RTX 3060 Laptop GPU with $6$GB of VRAM. A typical training epoch took $0.25$--$1$ seconds, giving a total training time of at most about five minutes. Pseudospectral grids were evaluated in batches of $512$; a batch typically took $1$--$2$ seconds, for a total time of about $30$ seconds to one minute. Trajectory generation was almost instant.

\subsection{Pendulum}

\begin{figure}[t]
\centering
\includegraphics[width=0.33\textwidth]{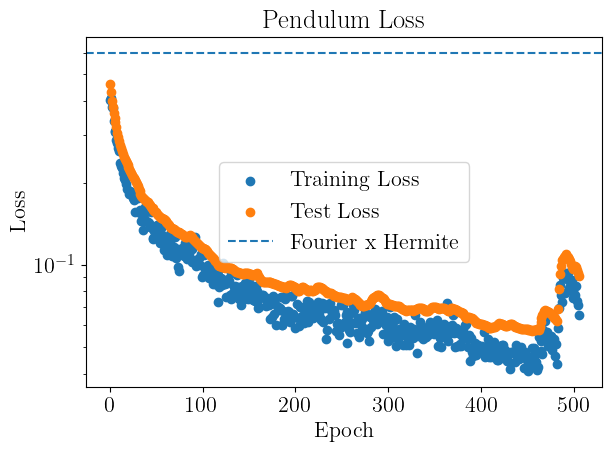}\hfill
\includegraphics[width=0.33\textwidth]{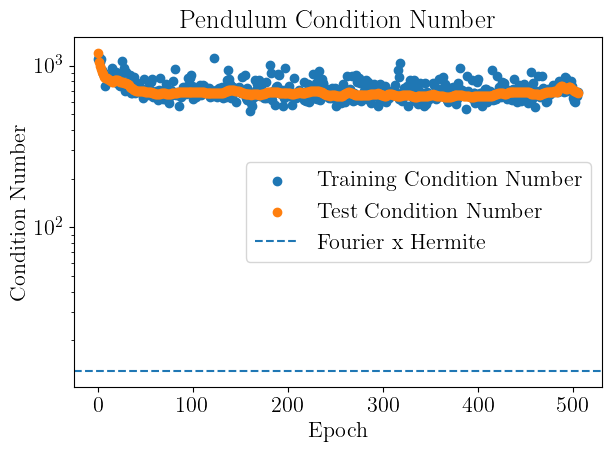}\hfill
\includegraphics[width=0.33\textwidth]{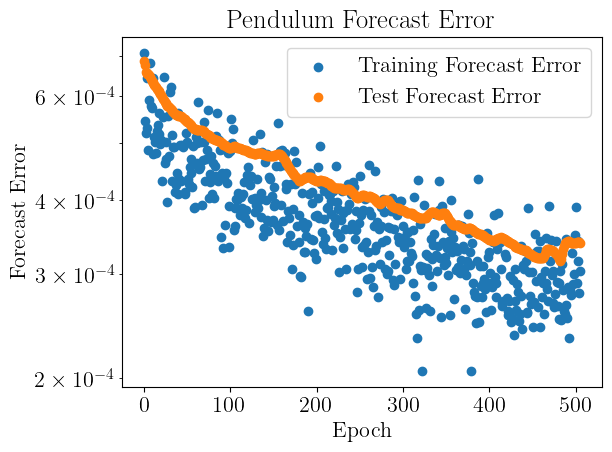}
\caption{Pendulum training diagnostics. Left: Training and test loss over $505$ epochs, with the Fourier--Hermite tensor dictionary shown as a horizontal reference line. Middle: Condition number over training and test. Right: Forecast error over training and test.}
\label{fig:pendulum_loss_cond}
\includegraphics[width=0.33\textwidth]{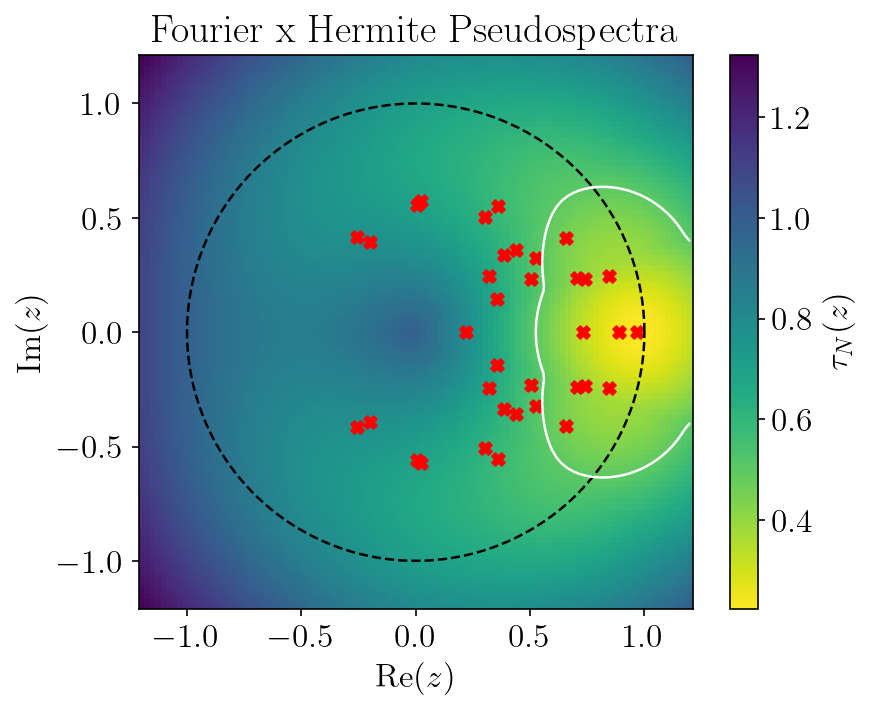}\hfill
\includegraphics[width=0.33\textwidth]{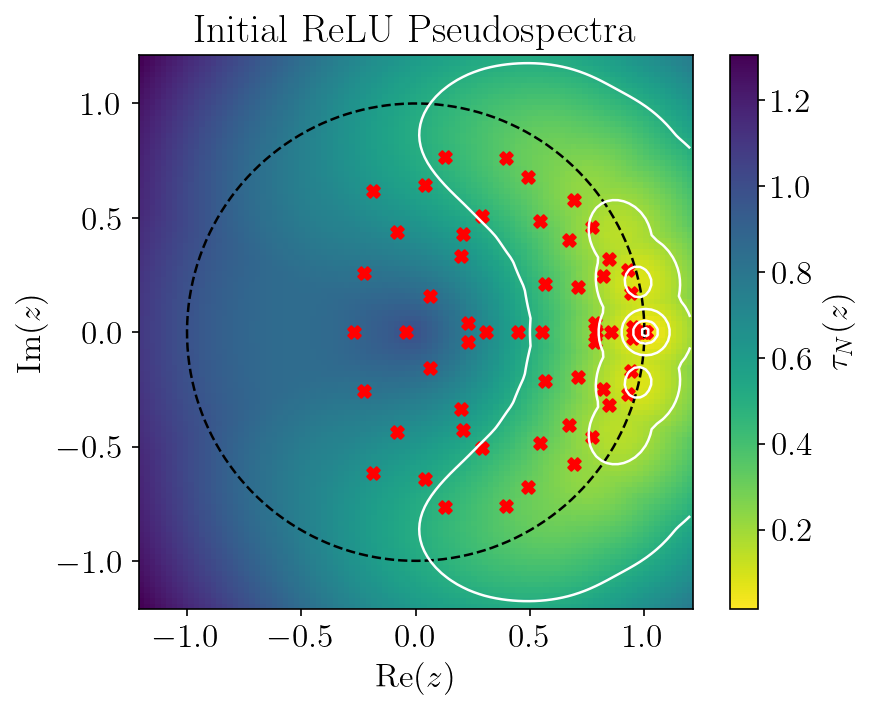}\hfill
\includegraphics[width=0.33\textwidth]{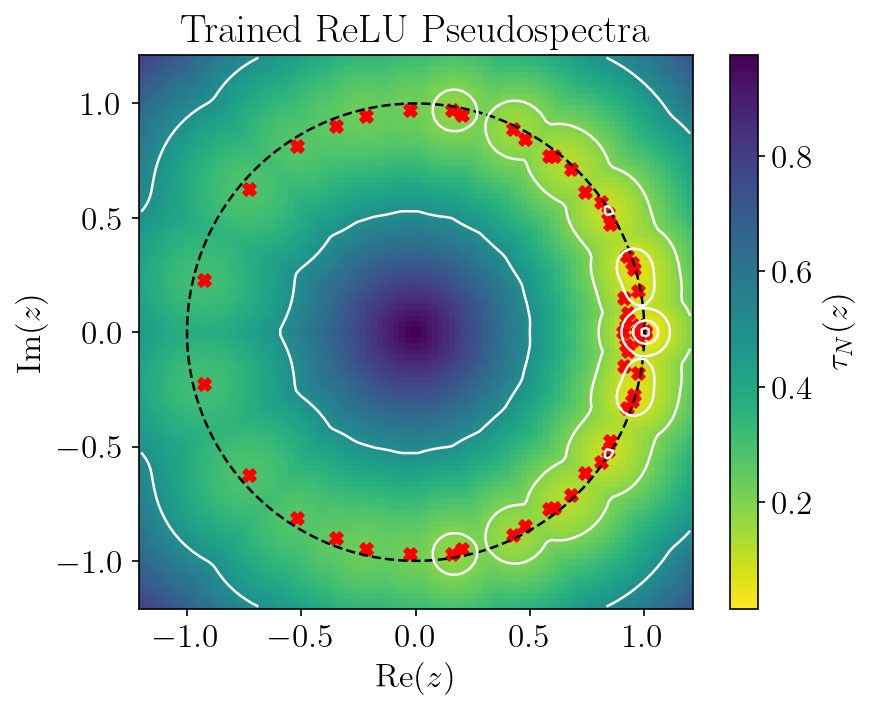}
\caption{Pendulum residual pseudospectra. Left: Fourier--Hermite dictionary $\{(\theta,p)\mapsto\varphi_j(\theta)h_k(p):0\leq j\leq 9,\ 0\leq k\leq 5\}$. Middle and right: ReLU-based dictionary before and after $505$ training epochs.}
\label{fig:pendulum}
\end{figure}

Our first example is the nonlinear pendulum. The state is $(\theta,p)\in\mathbb T\times\mathbb R$, where $\theta$ is the angle from the vertical and $p=\dot\theta$ is the angular velocity. With unit length and no damping, the equations are
$$
    \dot\theta=p,\qquad \dot p=-g\sin\theta,\qquad g=9.81.
$$
The continuous Hamiltonian flow preserves phase-space measure, and the associated Koopman operator is unitary. We therefore expect the computed spectrum to concentrate near $\mathbb T$. We use SciPy's \verb+solve_ivp+ to sample the solution with time step $\Delta t=0.1$. Trajectory data are generated from $140$ initial conditions $(\theta_j,p_j)\in[-\pi,\pi]\times[-1,1]$, with $\theta_j$ and $p_j$ sampled independently and uniformly from their respective intervals. Each initial condition is evolved for $126$ time steps, giving $17{,}640$ states in total. Training is run for at most $1000$ epochs, with early stopping after $50$ epochs without improvement in the training loss.

We benchmark the learned dictionary against the size-$60$ Fourier--Hermite dictionary
$$
    \bigl\{(\theta,p)\mapsto \varphi_j(\theta)h_k(p):0\leq j\leq 9,\ 0\leq k\leq 5\bigr\},
$$
where $(\varphi_j)_{j=0}^9$ are the first ten real Fourier modes in the angle coordinate and $h_k$ denotes the $k$th Hermite function. This choice reflects the geometry of the problem: eigenfunctions should be periodic in $\theta$, while the velocity coordinate is naturally treated on an unbounded domain.

\cref{fig:pendulum_loss_cond} shows the decrease in loss, condition number, and forecast error due to training. Over $10$ runs, the trained dictionary achieved loss $(6\pm0.5)\times10^{-2}$, condition number $512\pm68.7$, and forecast error $(5.03\pm0.57)\times10^{-4}$. By contrast, the Fourier--Hermite dictionary had loss around $0.6$ and forecast error around $0.62$.

\cref{fig:pendulum} compares the ResDMD residual pseudospectra and EDMD eigenvalues. The yellow regions are not merely visual features of the finite matrix: by the residual inclusion above, they identify spectral parameters with small Koopman residuals, and for the unitary pendulum Koopman operator this gives a direct distance-to-spectrum certificate. With the randomly initialized neural-network dictionary, several reliable modes appear near $z=1$, but there is significant spectral pollution closer to the centre of the disk and little spectral coverage on the left semicircle. After training, spectral coverage extends around the circle, with a substantial reduction in pollution.

\subsection{Undamped Harmonic Oscillator}

The undamped harmonic oscillator is the simplest conservative linear system. Let $x$ denote displacement from equilibrium and let $v$ denote velocity. With unit mass, the equations are
$$
    \dot x=v,\qquad \dot v=-kx,\qquad k=1.
$$
The exact continuous flow preserves phase-space measure, and the Koopman operator is unitary. We therefore regard computed eigenvalues far from $\mathbb T$ as spurious.

The system admits the exact solution
$$
    \begin{pmatrix}
        x(t)\\
        v(t)
    \end{pmatrix}
    =
    \begin{pmatrix}
        \cos t & \sin t\\
        -\sin t & \cos t
    \end{pmatrix}
    \begin{pmatrix}
        x_0\\
        v_0
    \end{pmatrix},
$$
where $x(0)=x_0$ and $v(0)=v_0$. We take $\Delta t=0.5$, dictionary size $N=36$, and hidden dimension $d_h=32$. Trajectory data are generated from $300$ initial conditions sampled uniformly from $[-1,1]^2$. Each initial condition is evolved for $60$ time steps, giving $18{,}000$ states in total. We benchmark the learned dictionary against the size-$36$ Hermite tensor dictionary
$$
    \bigl\{(x,v)\mapsto h_n(x)h_m(v):0\leq n\leq 5,\ 0\leq m\leq 5\bigr\},
$$
where $h_n$ denotes the $n$th Hermite function.

\begin{figure}[t]
\centering
\includegraphics[width=0.33\textwidth]{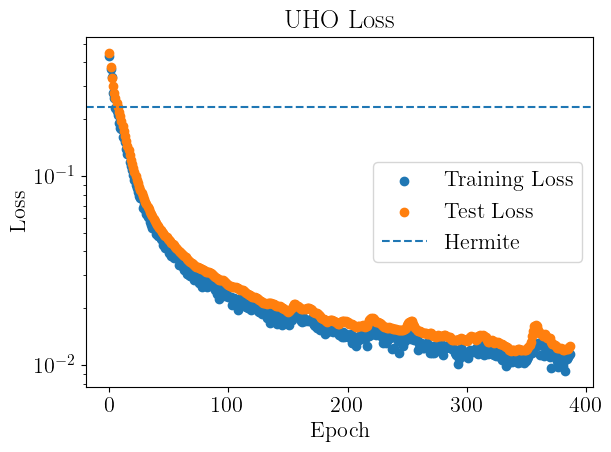}\hfill
\includegraphics[width=0.33\textwidth]{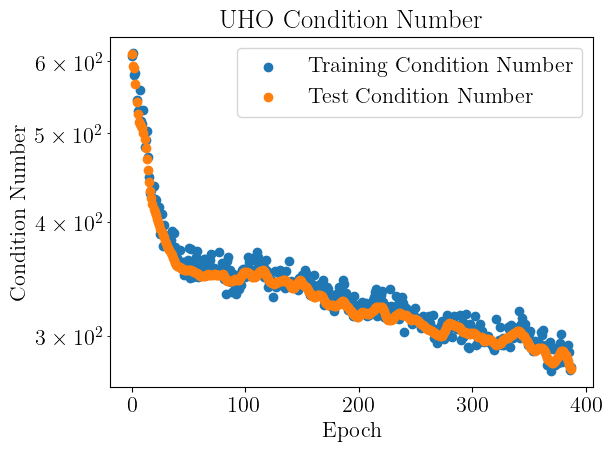}\hfill
\includegraphics[width=0.33\textwidth]{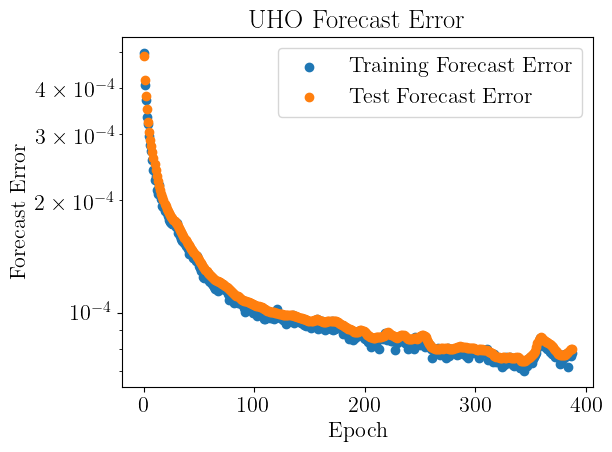}
\caption{Undamped harmonic oscillator diagnostics. Left: Training and test loss over $387$ epochs, with the Hermite dictionary shown as a horizontal reference line. Middle: Condition number over training and test. Right: Forecast error over training and test.}
\label{fig:undamped_loss_final}
\vspace{2mm}
\includegraphics[width=0.33\textwidth]{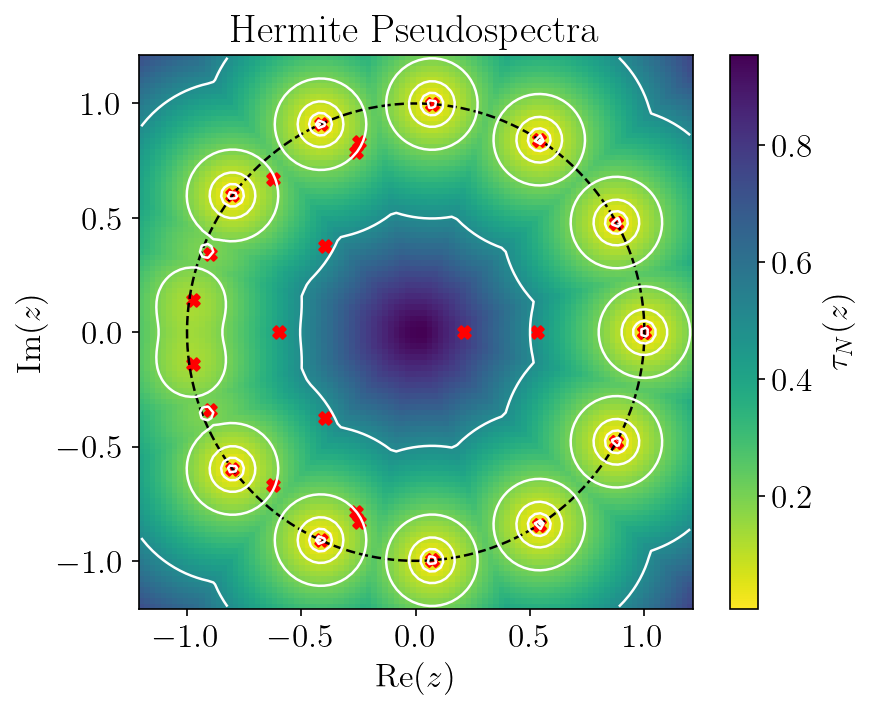}\hfill
\includegraphics[width=0.33\textwidth]{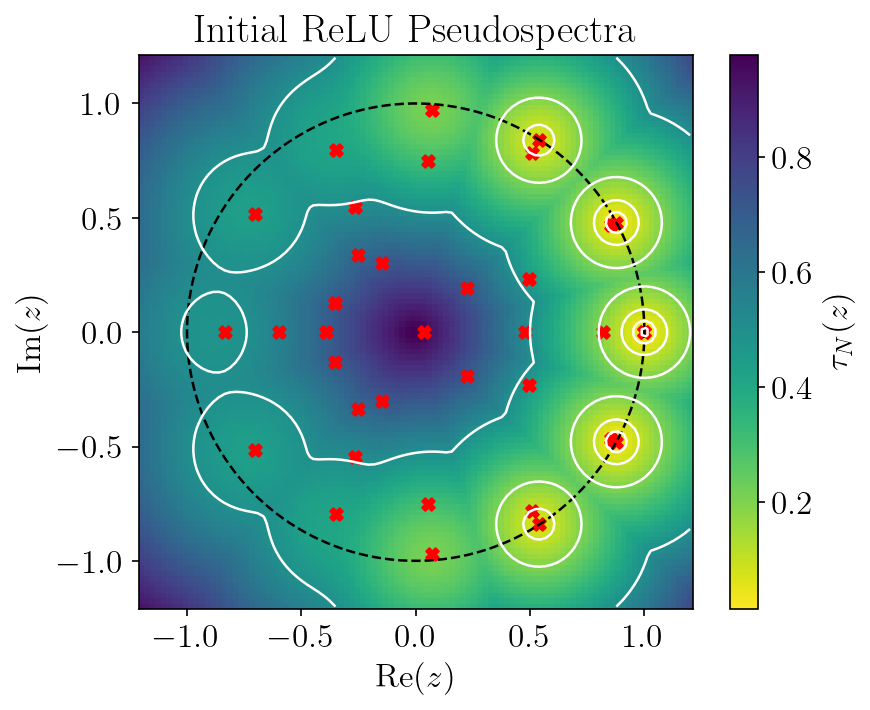}\hfill
\includegraphics[width=0.33\textwidth]{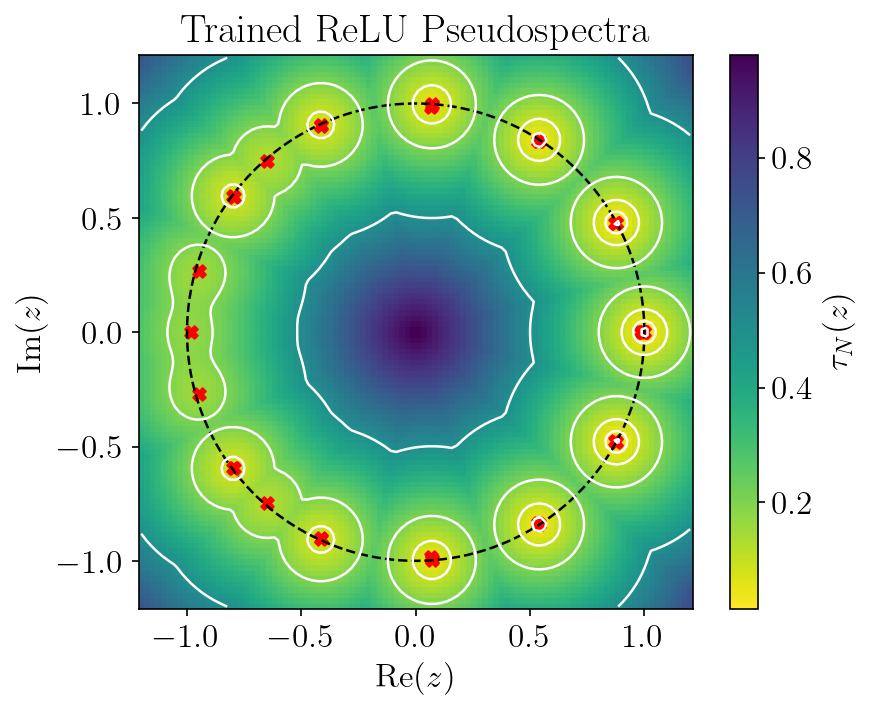}
\caption{Undamped harmonic oscillator residual pseudospectra. Left: Hermite dictionary. Middle and right: ReLU-based dictionary before and after $387$ training epochs.}
\label{fig:undamped_final}
\end{figure}

\cref{fig:undamped_loss_final} shows the decrease in loss, condition number, and forecast error due to training. Over $10$ runs, the trained dictionary achieved test loss $(1.04\pm0.01)\times10^{-2}$, condition number $258\pm40.4$, and test forecast error $(8.47\pm1.09)\times10^{-5}$. Its condition number improved meaningfully during training, from $634\pm54.6$ to $258\pm40.4$. By contrast, the Hermite dictionary had loss $0.232\pm0.003$ and condition number on the order of $10^6$. Nonetheless, it provided a reasonable approximation to the harmonic oscillator, with forecast error around $5\times10^{-3}$.

\cref{fig:undamped_final} shows that, before training, the neural-network dictionary already gives reasonable spectral inclusion on the right semicircle, with seven visible modes. Training removes the spurious eigenvalues away from $\mathbb T$ and resolves additional modes on the left semicircle, although these have slightly larger residuals than the modes on the right. The Hermite dictionary gives comparable spectral inclusion, but introduces several spurious eigenvalues inside the unit circle.

\subsection{Duffing Oscillator}

The Duffing oscillator is the dissipative system
$$
    \dot x=y,\qquad \dot y=-0.3y+x-x^3.
$$
We again use SciPy's \verb+solve_ivp+ to sample the continuous dynamics at discrete time steps. We use time step $\Delta t=0.05$ and dictionary size $N=100$. We sample $1500$ initial conditions from the box $[-3,3]^2$ and evolve each for $40$ time steps, producing $60{,}000$ snapshots in total. As a benchmark, we use the size-$100$ Chebyshev tensor dictionary
$$
    \bigl\{(x,y)\mapsto T_n(L_1^{-1}x)T_m(L_2^{-1}y):0\leq n\leq 9,\ 0\leq m\leq 9\bigr\},
$$
where $T_n$ denotes the $n$th Chebyshev polynomial, and the constants $L_1$ and $L_2$ are chosen so that the rescaled snapshots $(x^{(m)}/L_1,y^{(m)}/L_2)$ lie in $[-1,1]^2$. 

\begin{figure}[t]
\centering
\includegraphics[width=0.33\textwidth]{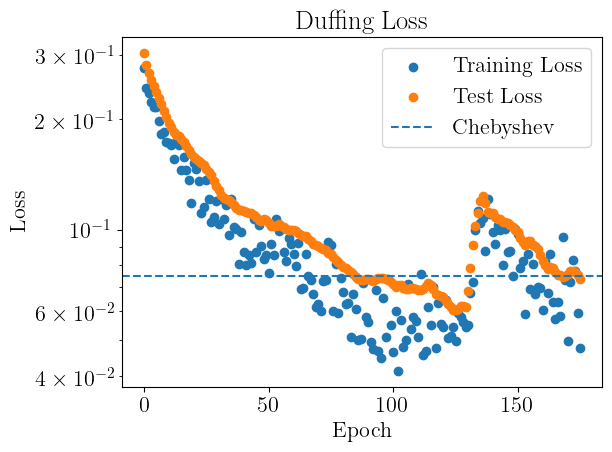}\hfill
\includegraphics[width=0.33\textwidth]{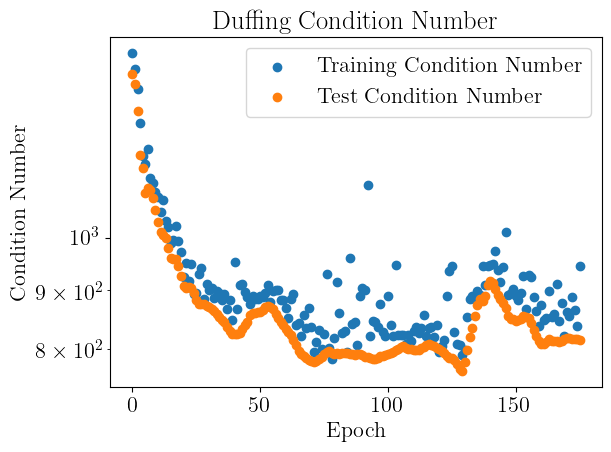}\hfill
\includegraphics[width=0.33\textwidth]{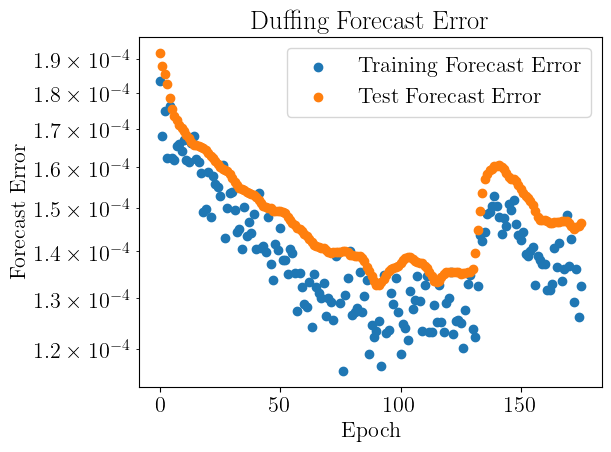}
\caption{Duffing oscillator diagnostics. Left: Training and test loss over $175$ epochs, with the benchmark dictionary shown as a horizontal reference line. Middle: Condition number over training and test. Right: Forecast error over training and test.}
\label{fig:duffing_loss_final}
\vspace{2mm}
\includegraphics[width=0.33\textwidth]{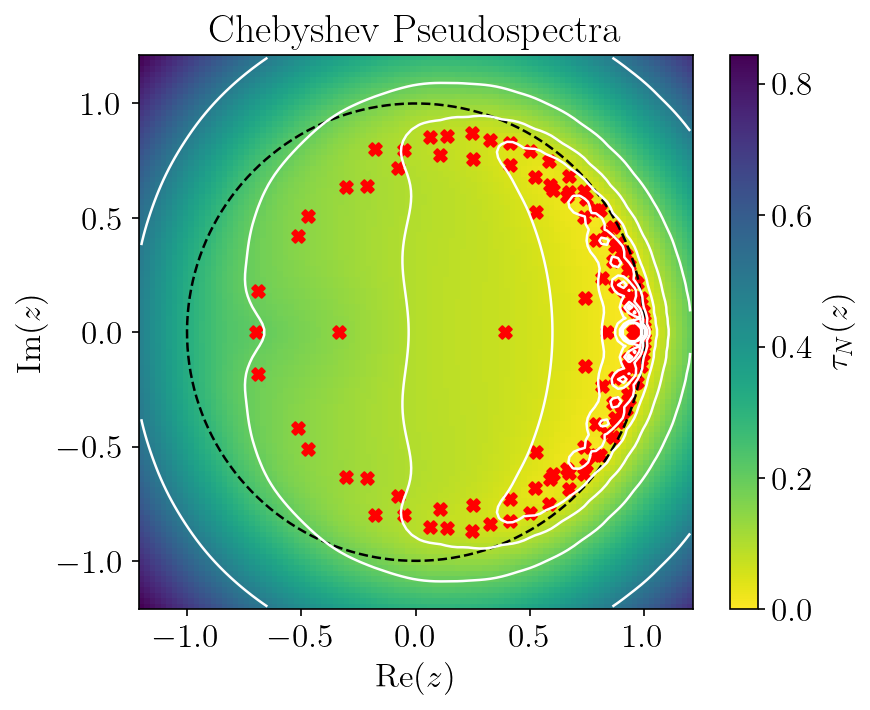}\hfill
\includegraphics[width=0.33\textwidth]{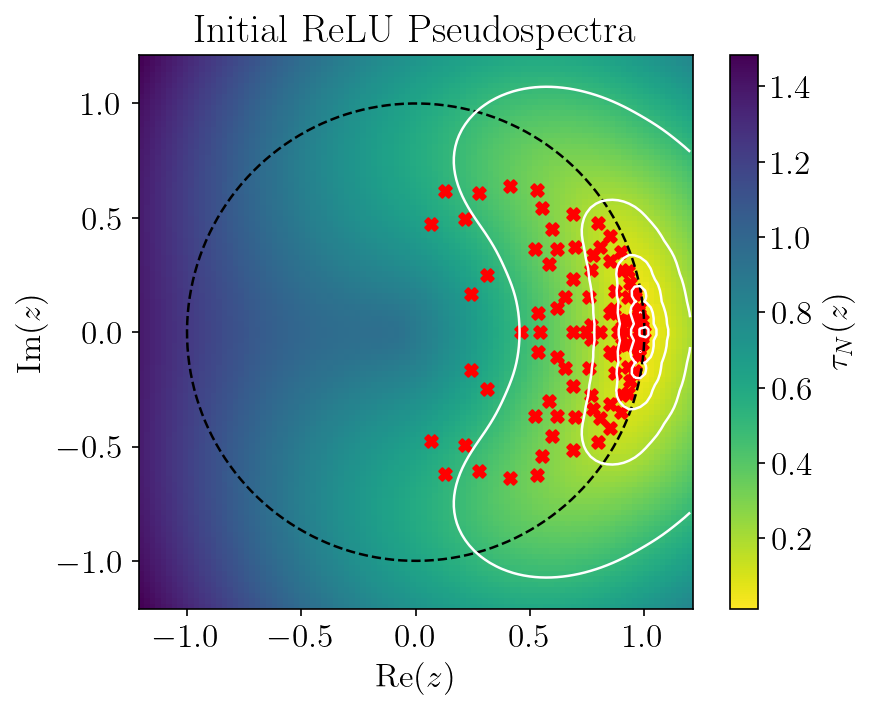}\hfill
\includegraphics[width=0.33\textwidth]{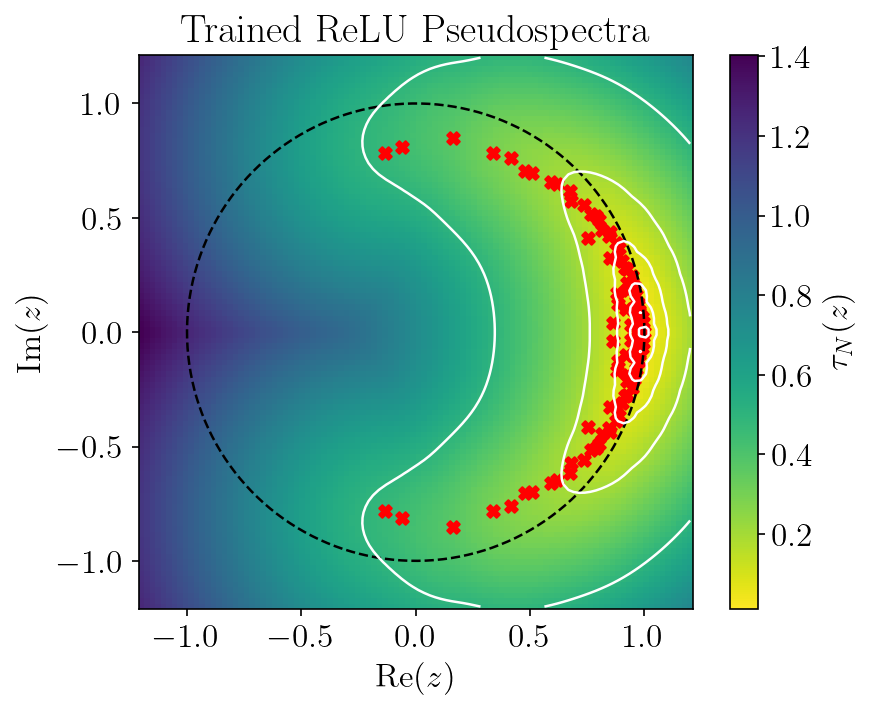}
\caption{Duffing oscillator residual pseudospectra. Left: Benchmark dictionary. Middle and right: ReLU-based dictionary before and after $175$ training epochs.}
\label{fig:duffing_final}
\end{figure}

\cref{fig:duffing_loss_final} shows the decrease in loss, condition number, and forecast error due to training. Over $10$ runs, the trained dictionary achieved loss $(3.03\pm0.47)\times10^{-2}$, condition number $717\pm20.6$, and test forecast error $(1.24\pm0.12)\times10^{-4}$. The Chebyshev dictionary had loss around $0.08$ and forecast error $(2.1\pm0.2)\times10^{-3}$. Its condition number was substantially larger, on the order of $10^5$.

A representative set of pseudospectral plots is shown in \cref{fig:duffing_final}. The Chebyshev dictionary exhibits significant spectral pollution throughout the disk. The untrained neural-network dictionary also shows substantial pollution: many eigenvalues cluster near $z=1$, while others appear in the interior of the disk. After training, the low-residual regions are more sharply resolved and the spurious interior eigenvalues are substantially reduced.

\section{Sea-Surface Temperature}
We conclude with an example based on real-world climate data, to test whether the improvements above persist for noisy observations not generated from an exact model. We consider monthly sea-surface temperature data.

We use the NOAA OISST v2.1\footnote{\url{https://www.ncei.noaa.gov/data/sea-surface-temperature-optimum-interpolation/v2.1/access/avhrr/}} dataset, which gives a complete map of ocean surface temperature at $0.25^\circ\times0.25^\circ$ spatial resolution. We restrict to the tropical Pacific,
$$
    120^\circ\mathrm E\leq\lambda\leq 80^\circ\mathrm W,\qquad 20^\circ\mathrm S\leq\phi\leq 20^\circ\mathrm N,
$$
where $\lambda$ is longitude and $\phi$ is latitude. We use monthly data from September 1981 through May 2026, with provisional values used for the final two weeks.

Let $S(t,\xi)$ denote the mean sea-surface temperature in month $t$ at grid point $\xi$. The month $t=1$ is September 1981, and $t=537$ is May 2026. For a month $t$, let $m(t)\in\{1,\ldots,12\}$ denote its calendar month, so that $m(1)=9$ and $m(537)=5$. Let $\overline S_{m(t)}(\xi)$ be the climatological average at grid point $\xi$ for calendar month $m(t)$. We form the anomalies
$$
    A(t,\xi)=S(t,\xi)-\overline S_{m(t)}(\xi).
$$

A naive approach would form a roughly $95{,}000$-dimensional vector from the values $A(t,\xi)$ over the full grid. This is computationally expensive, and the $537$ monthly snapshots are too few to reliably resolve such a high-dimensional system. Instead, we perform principal component analysis (PCA). We apply the standard latitude weighting $[\cos(\phi_\xi)]^{1/2}$, where $\phi_\xi$ is the latitude of grid point $\xi$. On a longitude-latitude grid, cells with equal angular resolution have physical area proportional to $\cos(\phi_\xi)$; the weighting therefore accounts for unequal cell areas.

PCA decomposes the anomaly field as a linear combination of empirical orthogonal functions,
$$
    A(t,\xi)\approx\sum_{j=1}^r a_j(t)e_j(\xi).
$$
The first $30$ principal components explain around $90\%$ of the variance, with the first $10$ explaining around $80\%$ and the first $5$ explaining around $70\%$. We therefore run experiments with $r=5$ and $r=10$. The reduced states are
$$
    x_t=(a_1(t),\ldots,a_r(t))\in\mathbb R^r,\qquad r=5,10,
$$
and the snapshot pairs are $(x_t,x_{t+1})$, corresponding to one-month transitions.

The ground truth for spectral accuracy is less clear than in the preceding benchmark systems, which are theoretically well understood. We therefore place more emphasis on forecast error when assessing the relevance of the learned dictionary. To prevent data leakage, we use a sequential rather than random split, with the first $70\%$ of the data used for training and the remaining $30\%$ used for testing. We compare the learned dictionaries with quadratic polynomial dictionaries containing $1$, $x_i$, and $x_ix_j$ for $1\leq i,j\leq r$. As before, trained results are averaged over $15$ runs, with standard deviations reported.

The number of snapshots is small compared with the synthetic examples, which typically had on the order of $10{,}000$ snapshots. It is therefore important to choose dictionaries that are expressive enough to capture the dominant dynamics, but not so expressive that they overfit. Indeed, even with modest dictionary sizes and architectures, the test loss may remain roughly flat while the training loss decreases. We therefore use relatively small trained dictionaries. Since the train-test split is deterministic, variation between runs comes only from random initialization; statistics for the fixed polynomial dictionaries are deterministic.

\begin{figure}[t]
\centering
\includegraphics[width=0.45\textwidth]{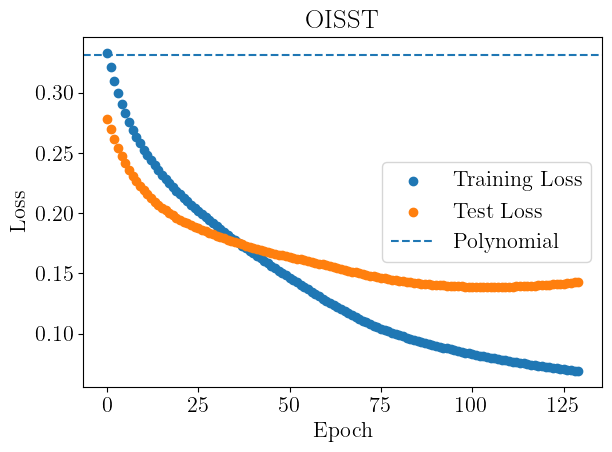}\hfill
\includegraphics[width=0.45\textwidth]{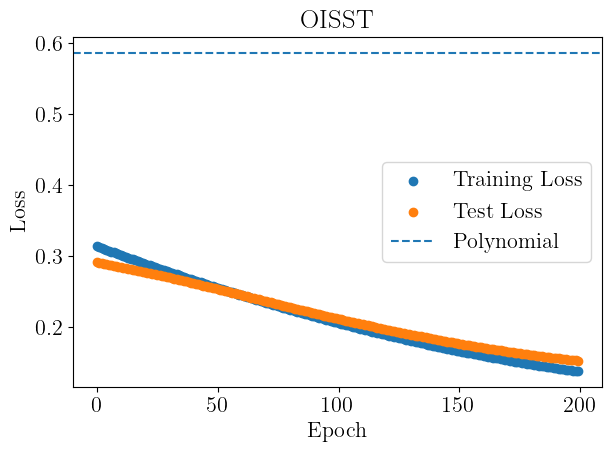}
\caption{SST loss curves for typical runs with $r=5$ and $r=10$. The polynomial dictionary loss is shown as a horizontal reference line.}
\label{fig:climate_loss}
\end{figure}

\begin{figure}[t]
\centering
\includegraphics[width=0.45\textwidth]{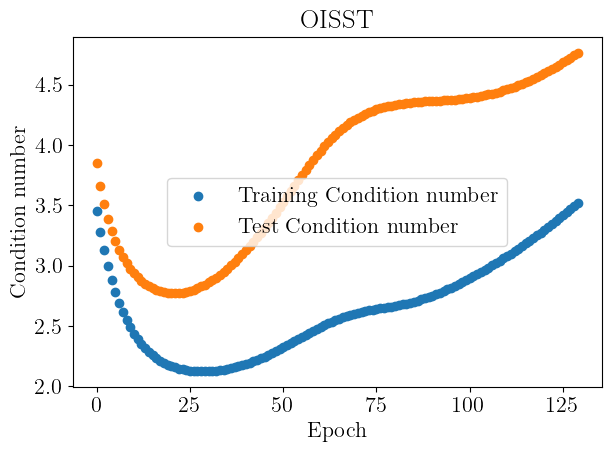}\hfill
\includegraphics[width=0.45\textwidth]{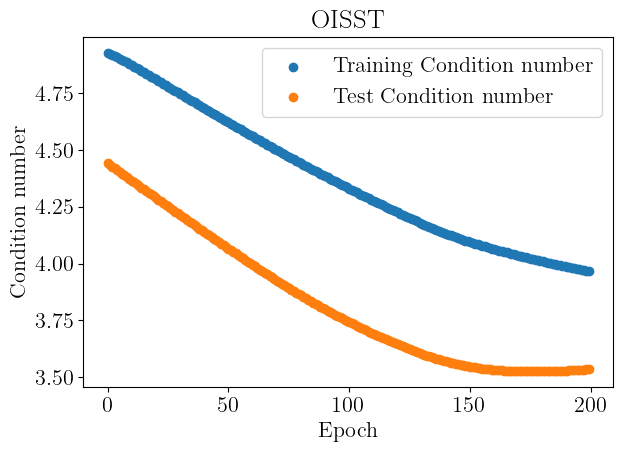}
\caption{SST condition-number curves for typical runs with $r=5$ and $r=10$.}
\label{fig:climate_cond_num}
\end{figure}

\begin{figure}[t]
\centering
\includegraphics[width=0.45\textwidth]{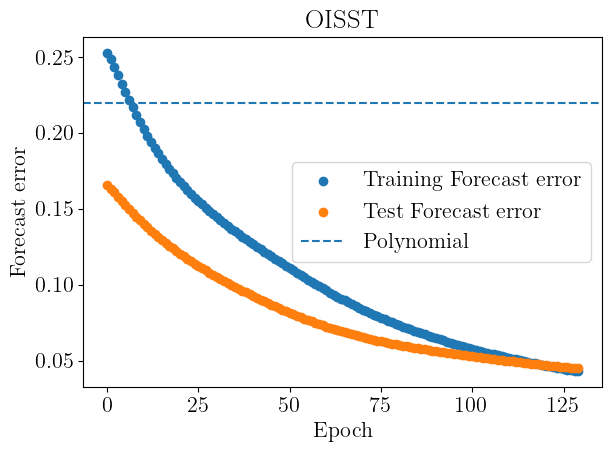}\hfill
\includegraphics[width=0.45\textwidth]{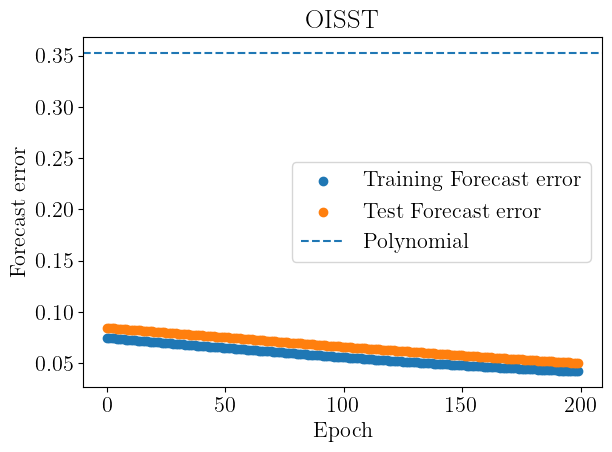}
\caption{SST forecast-error curves for typical runs with $r=5$ and $r=10$. The polynomial dictionary forecast error is shown as a horizontal reference line.}
\label{fig:climate_forecast_error}
\end{figure}

With $r=5$, using dictionary size $4$, one hidden layer, and hidden dimension $64$, the trained dictionary achieves loss $0.1584\pm0.0232$, condition number $5.7298\pm1.3507$, and forecast error $0.0433\pm0.0090$. The small condition number is partly due to the small dictionary size. The quadratic dictionary achieves loss $0.3111$, condition number $7.3202$, and forecast error $0.2194$. With $r=10$, using dictionary size $3$, one hidden layer, and hidden dimension $16$, the trained dictionary achieves loss $0.2136\pm0.0501$, condition number $5.4797\pm2.0335$, and forecast error $0.0569\pm0.0277$. The quadratic dictionary achieves loss $0.5866$, condition number $13.5065$, and forecast error $0.3529$. The loss, condition number, and forecast error curves are shown in \cref{fig:climate_loss,fig:climate_cond_num,fig:climate_forecast_error}.

\begin{figure}[t]
\centering
\includegraphics[width=1\textwidth]{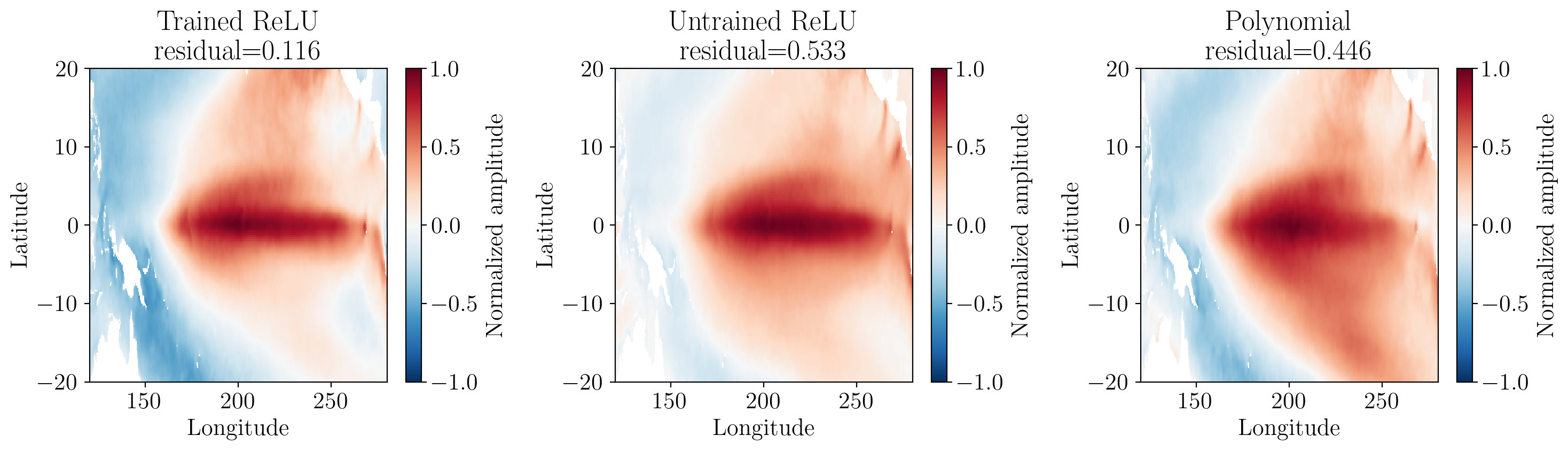}
\caption{Koopman mode with strongest amplitude in the Ni\~{n}o 3.4 region. Training substantially reduces the associated residual.}
\label{fig:climate_nino}
\end{figure}

The Ni\~{n}o 3.4 region spans $120^\circ\mathrm W$ to $170^\circ\mathrm W$ longitude and $5^\circ\mathrm S$ to $5^\circ\mathrm N$ latitude. Sea-surface temperature in this region is used to monitor El Ni\~{n}o and La Ni\~{n}a events. We select the Koopman modes with the strongest amplitude in this region and compare their residuals before and after training, as well as against the polynomial dictionary. For $r=5$, training gives a significant reduction in residual, as shown in \cref{fig:climate_nino}. Amplitudes are normalized to lie in $[-1,1]$.

\section{Conclusion}

Koopman approximation is not just a problem of prediction. It is a problem of spectral meaning. A finite EDMD matrix always has eigenvalues, but without an operator-level test they may be little more than artifacts of the chosen dictionary. ResDMD supplies the missing test: residuals for the Koopman operator, not just for its finite projection.

This paper has used that test as a learning signal. Neural networks provide flexible dictionaries, but flexibility alone is dangerous. Residual minimization drives the dictionary toward spectrally relevant observables; conditioning keeps those observables numerically distinct. The two requirements belong together. Without the first, the spectrum is not trustworthy. Without the second, the computation is not stable.

The experiments support this view. For conservative systems, training improves residual pseudospectral inclusion and suppresses spurious eigenvalues away from the unit circle. For the dissipative Duffing oscillator, it sharpens low-residual spectral structure while reducing pollution in the disk. For sea-surface temperature data, where there is no exact model to appeal to, the same principle gives lower residuals and better forecasts than a standard polynomial dictionary.

The conclusion is simple. Learned dictionaries should not be judged only by how well they predict the next snapshot. They should be judged by whether they produce Koopman objects that can be trusted. Residuals give the spectral target; conditioning makes the computation viable. Together they turn neural Koopman learning from a flexible fitting procedure into a spectrally disciplined numerical method.

\section*{Code availability}
The implementation used to generate the dictionaries, train the neural networks and reproduce the experiments will be made available at \url{https://github.com/GeorgeCoote/pyresdmd}.

\bibliographystyle{plain}
\bibliography{references}

%USE THE BELOW OPTIONS IN CASE YOU NEED AUTHOR YEAR FORMAT.
%\bibliographystyle{abbrvnat}
%\bibliography{references}

\end{document}